\documentclass[twoside,12pt]{article}
\usepackage{amssymb}
\usepackage{amsfonts}
\usepackage{amsmath}

\input{tcilatex}

\begin{document}

\title{Cohomology associated to a Poisson structure on Weil bundles \\
\medskip }
\author{Vann Borhen NKOU$^{1,2}$; Basile Guy Richard BOSSOTO$^{1}$ \and %
\medskip \\
$^{1}$Marien NGOUABI University, Faculty of Sciences\\
Department of Mathematics\\
BP: 69 - Brazzaville (Congo)\\
E-mail: bossotob@yahoo.fr\\
\medskip \\
$^{2}$Abomey Calavi University, IMSP\\
BP: 13, Porto-novo, Benin.\\
E-mail: vannborhen@yahoo.fr}
\date{}
\maketitle

\begin{abstract}
Let $M$ be a paracompact smooth manifold of dimension $n,$ $A$ a Weil
algebra and $M^{A}$ the Weil bundle associated. We define and describe the
notion of $\widetilde{d}$-Poisson cohomology and of $\widetilde{d_{A}}$%
-Poisson cohomology on $M^{A}$.
\end{abstract}



\section{Introduction}

A local algebra in the sens of Andr\'{e} Weil or simply a Weil algebra is a
real, unitary, commutative algebra of finite dimension with a unique maximal
ideal of codimension $1$ on $\mathbb{R}$ \cite{wei}.\newline
Let $A$ be a Weil algebra and $\mathfrak{m}$ be its maximal ideal. We have 
\begin{equation*}
A=\mathbb{R}\oplus \mathfrak{m}.
\end{equation*}%
The first projection 
\begin{equation*}
A=\mathbb{R}\oplus \mathfrak{m}\longrightarrow \mathbb{R}
\end{equation*}%
is a homomorphism of algebra which is surjective, called augmentation and
the unique none-zero integer $k\in \mathbb{N}$ such that $\mathfrak{m}%
^{k}\neq (0)$ and $\mathfrak{m}^{k+1}=(0)$ is the height of $A$.\newline

If $M$ is a smooth manifold, $C^{\infty }(M)$ the algebra of differentiable
functions on $M$ and $A$ a Weil algebra of maximal ideal $\mathfrak{m}$, an
infinitely near point to $x\in M$ of kind $A$ is a homomorphism of algebras 
\begin{equation*}
\xi :C^{\infty }(M)\longrightarrow A
\end{equation*}%
such that $[\xi (f)-f(x)]\in \mathfrak{m}$ for any $f\in C^{\infty }(M)$.%
\newline

We denote $M_{x}^{A}$ the set of all infinitely near points to $x\in M$ of
kind $A$ and 
\begin{equation*}
M^{A}=\bigcup\limits_{x\in M}M_{x}^{A}\text{.}
\end{equation*}

The set $M^{A}$ is a smooth manifold of dimension $\dim M\times \dim A$
called manifold of infinitely near points of kind $A$\cite{mor}.\newline

When both $M$ and $N$ are smooth manifolds and when 
\begin{equation}
h:M\longrightarrow N  \notag
\end{equation}%
is a differentiable application, then the application 
\begin{equation}
h^{A}:M^{A}\longrightarrow N^{A},\xi \longmapsto h^{A}(\xi ),  \notag
\end{equation}%
such that, for any $g\in C^{\infty }(N)$, 
\begin{equation}
\left[ h^{A}(\xi )\right] (g)=\xi (g\circ h)  \notag
\end{equation}%
is also differentiable. When $h$ is a diff{e}omorphism, it is the same for $%
h^{A}$.\newline

Moreover, if $\varphi :A\longrightarrow B$ is a homomorphism of Weil
algebras, for any smooth manifold $M$, the application 
\begin{equation}
\varphi _{M}:M^{A}\longrightarrow M^{B},\xi \longmapsto \varphi \circ \xi 
\notag
\end{equation}%
is differentiable. In particular, the augmentation 
\begin{equation}
A\longrightarrow \mathbb{R}  \notag
\end{equation}%
defines for any smooth manifold $M$, the projection 
\begin{equation}
\pi _{M}:M^{A}\longrightarrow M,  \notag
\end{equation}%
which assigns every infinitely near point to $x\in M$ to its origin $x$.
Thus ($M^{A},\pi _{M},M$) defines the bundle of infinitely near points or
simply weil bundle \cite{oka},\cite{kms},\cite{wei}.

If $(U,\varphi )$ is a local chart of $M$ with coordinate functions $%
(x_{1},x_{2},...,x_{n})$, the application 
\begin{equation}
U^{A}\longrightarrow A^{n},\xi \longmapsto (\xi (x_{1}),\xi (x_{2}),...,\xi
(x_{n})),  \notag
\end{equation}%
is a bijection from $U^{A}$ into an open of $A^{n}$. The manifold $M^{A}$ is
a smooth manifold modeled over $A^{n}$, that is to say an $A$-manifold of
dimension $n$ \cite{bok},\cite{shu}.\newline

The set, $C^{\infty }(M^{A},A)$ of differentiable functions on $M^{A}$ with
value in $A$ is a commutative, unitary algebra over $A$.When one
identitifies $\mathbb{R}^{A}$ with $A$, for $f\in C^{\infty }(M)$, the
application%
\begin{equation}
f^{A}:M^{A}\longrightarrow A,\xi \longmapsto \xi (f)  \notag
\end{equation}%
is differentiable. Moreover the application%
\begin{equation}
C^{\infty }(M)\longrightarrow C^{\infty }(M^{A},A),f\longmapsto f^{A}, 
\notag
\end{equation}%
is an injective homomorphism of algebras and we have: 
\begin{equation}
(f+g)^{A}=f^{A}+g^{A};(\lambda \cdot f)^{A}=\lambda \cdot f^{A};(f\cdot
g)^{A}=f^{A}\cdot g^{A}  \notag
\end{equation}%
for $\lambda \in \mathbb{R}$, $f$ and $g$ belonging to $C^{\infty }(M)$.%
\newline

We denote $\mathfrak{X}(M^{A})$, the set of all vector fields on $M^{A}.$
According to \cite{bok}, We have the following equivalent assertions:

\begin{enumerate}
\item $X:C^{\infty }(M^{A})\longrightarrow C^{\infty }(M^{A})$ is a vector
field on $M^{A}$;

\item $X:C^{\infty }(M)\longrightarrow C^{\infty }(M^{A},A)$ is a linear
application which verifies%
\begin{equation*}
X(fg)=X(f)\cdot g^{A}+f^{A}\cdot X(g)
\end{equation*}%
for any $f,g\in C^{\infty }(M)$ i.e is a derivation of $C^{\infty }(M)$ into 
$C^{\infty }(M^{A},A)$ with respect to the module structure 
\begin{equation*}
C^{\infty }(M^{A},A)\times C^{\infty }(M)\longrightarrow C^{\infty
}(M^{A},A),(\varphi ,f)\longmapsto \varphi \cdot f^{A}\text{.}
\end{equation*}%
Thus, the set $\mathfrak{X}(M^{A})$ of all vector fields on $M^{A}$ is a $%
C^{\infty }(M^{A},A)$-module.\newline
\end{enumerate}

\bigskip When%
\begin{equation}
\theta :C^{\infty }(M)\longrightarrow C^{\infty }(M)  \notag
\end{equation}%
is a vector field on $M$, the application%
\begin{equation}
\theta ^{A}:C^{\infty }(M)\longrightarrow C^{\infty }(M^{A},A),f\longmapsto 
\left[ \theta (f)\right] ^{A}\text{,}  \notag
\end{equation}%
is a vector field on $M^{A}$. The vector field $\theta ^{A}$ is the
prolongation to $M^{A}$ of the vector field $\theta $.\newline

\begin{theorem}
If $X$ is a vector field on $M^{A}$ considered as a derivation of $C^{\infty
}(M)$ into $C^{\infty}(M^{A},A)$, then there exists, an unique derivation 
\begin{equation*}
\widetilde{X}:C^{\infty}(M^{A},A)\longrightarrow C^{\infty}(M^{A},A) 
\end{equation*}
such that such that

\begin{enumerate}
\item $\widetilde{X}$ is $A$-linear;

\item $\widetilde{X}\left[ C^{\infty}(M^{A})\right] \subset C^{\infty}(M^{A})
$;

\item $\ \widetilde{X}(f^{A})=X(f)$ for any $f\in C^{\infty}(M)$.
\end{enumerate}
\end{theorem}

Thus, the application%
\begin{equation}
\lbrack ,]:\mathfrak{X}(M^{A})\times \mathfrak{X}(M^{A})\longrightarrow 
\mathfrak{X}(M^{A}),(X,Y)\longmapsto \widetilde{X}\circ Y-\widetilde{Y}\circ
X,  \notag
\end{equation}%
is $A$-bilinear and defines a structure of Lie algebra over $A$ on $%
\mathfrak{X}(M^{A})$\cite{bok}.\newline

The goal of this paper is to define and describe the notion of $\widetilde{d}
$-Poisson cohomology and of $\widetilde{d_{A}}$-Poisson cohomology.

\section{Poisson structure on Weil bundles}

In this section, $M$ is a Poisson manifold i.e there exists a bracket $\{,\}$
on $C^{\infty }(M)$ such that the pair $(C^{\infty }(M),\{,\})$ is a real
Lie algebra and for any $f\in C^{\infty }(M),$ the application%
\begin{equation}
ad(f):C^{\infty }(M)\longrightarrow C^{\infty }(M),g\longmapsto \{f,g\} 
\notag
\end{equation}%
is a derivation of commutative algebra i.e 
\begin{equation}
\{f,g\cdot h\}=\{f,g\}\cdot h+g\cdot \{f,h\}  \notag
\end{equation}%
$f,g,h\in C^{\infty }(M)$ \cite{Lic},\cite{iva}.

We denote%
\begin{equation}
C^{\infty }(M)\longrightarrow Der_{\mathbb{R}}[C^{\infty }(M)],f\longmapsto
ad(f)\text{,}  \notag
\end{equation}%
the adjoint representation and $d$ the operator of cohomology associated to
this representation. For any $p\in \mathbb{N}$, 
\begin{equation}
\Lambda _{Pois}^{p}(M)=C^{p}[C^{\infty }(M),C^{\infty }(M)]  \notag
\end{equation}%
denotes the $C^{\infty }(M)$-module of skew-symmetric multilinear forms of
degree $p$ from $C^{\infty }(M)$ into $C^{\infty }(M)$. We have 
\begin{equation}
\Lambda _{Pois}^{0}(M)=C^{\infty }(M).  \notag
\end{equation}

The $A$-algebra $C^{\infty }(M^{A},A)$ is a Poisson algebra over $A$ if
there exists a bracket $\{,\}$ on $C^{\infty }(M^{A},A)$ such that the pair $%
(C^{\infty }(M^{A},A),\{,\})$ is a Lie algebra over $A$ satisfying%
\begin{equation*}
\{\varphi _{1}\cdot \varphi _{2},\varphi _{3}\}=\{\varphi _{1},\varphi
_{3}\}\cdot \varphi _{2}+\varphi _{1}\cdot \{\varphi _{2},\varphi _{3}\}
\end{equation*}%
for any $\varphi _{1},\varphi _{2},\varphi _{3}\in C^{\infty }(M^{A},A)$ 
\cite{bok1},\cite{bos2}.$\newline
$When $M$ is a Poisson manifold with bracket $\{,\}$, for any $f,g\in
C^{\infty }(M)$,%
\begin{equation}
ad(fg)=ad(f)\cdot g+f\cdot ad(g)\text{.}  \notag
\end{equation}

For any $f\in C^{\infty}(M)$, let%
\begin{equation}
\lbrack ad(f)]^{A}:C^{\infty}(M)\longrightarrow
C^{\infty}(M^{A},A),g\longmapsto\{f,g\}^{A}\text{,}  \notag
\end{equation}
be the prolongation of the vector field $ad(f)$ and let%
\begin{equation}
\widetilde{\lbrack ad(f)]^{A}}:C^{\infty}(M^{A},A)\longrightarrow C^{\infty
}(M^{A},A)  \notag
\end{equation}
be the unique $A$-linear derivation such that%
\begin{equation*}
\widetilde{\lbrack ad(f)]^{A}}(g^{A})=[ad(f)]^{A}(g)=\{f,g\}^{A}
\end{equation*}
for any $g\in C^{\infty}(M)$.

\begin{theorem}
\cite{bok1} For $\varphi\in C^{\infty}(M^{A},A)$, the application 
\begin{equation*}
\tau_{\varphi}:C^{\infty}(M)\longrightarrow C^{\infty}(M^{A},A),f\longmapsto
-\widetilde{[ad(f)]^{A}}(\varphi) 
\end{equation*}
is a vector field on $M^{A}$.
\end{theorem}

We denote%
\begin{equation*}
\widetilde{\tau_{\varphi}}:C^{\infty}(M^{A},A)\longrightarrow
C^{\infty}(M^{A},A)
\end{equation*}
the unique $A$-linear derivation such that%
\begin{equation}
\widetilde{\tau_{\varphi}}(f^{A})=\tau_{\varphi}(f)  \notag
\end{equation}
for any $f\in C^{\infty}(M)$. We have for $f\in C^{\infty}(M)$, 
\begin{equation*}
\widetilde{\tau_{f^{A}}}=\widetilde{[ad(f)]^{A}}\text{,}
\end{equation*}
and for $\varphi,\psi\in C^{\infty}(M^{A},A)$ and for $a\in A,$ 
\begin{align}
\widetilde{\tau}_{\varphi+\psi} & =\widetilde{\tau}_{\varphi}+\widetilde {%
\tau}_{\psi};  \notag \\
\widetilde{\tau}_{a\cdot\varphi} & =a\cdot\widetilde{\tau}_{\varphi }; 
\notag \\
\widetilde{\tau}_{\varphi\cdot\psi} & =\varphi\cdot\widetilde{\tau}_{\psi
}+\psi\cdot\widetilde{\tau}_{\varphi}\text{.}  \notag
\end{align}
For any $\varphi,\psi\in C^{\infty}(M^{A},A)$, we let%
\begin{equation}
\{\varphi,\psi\}_{A}=\widetilde{\tau}_{\varphi}(\psi).  \notag
\end{equation}
In \cite{bok1} we show that this bracket defines a structure of $A$-Poisson
algebra on $C^{\infty}(M^{A},A)$.

\begin{theorem}
If $M$ is a Poisson manifold with bracket $\{,\}$, then $\{,\}_{A}$ is the
prolongation on $M^{A}$ of the structure of Poisson on $M$ defined by $\{,\}$%
.
\end{theorem}

\bigskip

\section{$\widetilde{d}$-Poisson cohomology}

\bigskip

\begin{proposition}
When $M$ is a Poisson manifold with bracket $\{,\}$, the map 
\begin{equation*}
C^{\infty}(M)\longrightarrow Der_{A}\left[ C^{\infty}(M^{A},A)\right]
,f\longmapsto\widetilde{[ad(f)]^{A}}
\end{equation*}
is a representation of $C^{\infty}(M)\ $into $C^{\infty}(M^{A},A)$.
\end{proposition}

We denote $\ \widetilde{d}$ the operator of cohomology associated to this
representation.\newline
For any $p\in \mathbb{N}$,%
\begin{equation}
\Lambda _{Pois}^{p}(M^{A},\sim )=\mathcal{C}^{p}[C^{\infty }(M),C^{\infty
}(M^{A},A)]  \notag
\end{equation}%
denotes the $C^{\infty }(M^{A},A)$-module of skew-symmetric multilinear
forms of degree $p$ from $C^{\infty }(M)$ into $C^{\infty }(M^{A},A)$. We
have 
\begin{equation}
\Lambda _{Pois}^{0}(M^{A},\sim )=C^{\infty }(M^{A},A).  \notag
\end{equation}%
We denote%
\begin{equation}
\Lambda _{Pois}(M^{A},\sim )=\bigoplus_{p=0}^{n}\Lambda
_{Pois}^{p}(M^{A},\sim )\text{.}  \notag
\end{equation}

Thus, for $\Omega \in \Lambda _{Pois}^{p}(M^{A},\sim )$ and $%
f_{1},...,f_{p+1}\in C^{\infty }(M),$ we have 
\begin{align*}
\widetilde{d}\Omega (f_{1},...,f_{p+1})& =\sum_{i=1}^{p+1}(-1)^{i}\widetilde{%
[ad(f_{i})]^{A}}[\Omega (f_{1},...,\widehat{f_{i}},...,f_{p+1})] \\
& +\sum_{1\leq i<j\leq p+1}(-1)^{i+j}\Omega (\{f_{i},f_{j}\},f_{1},...,%
\widehat{f_{i}},...,\widehat{f_{j}},...,f_{p+1})
\end{align*}%
where $\widehat{f_{i}}$ means that the term $f_{i}$ is omitted.\newline

When $\eta \in \Lambda _{Pois}^{p}(M)$, then%
\begin{equation}
\eta ^{A}:C^{\infty }(M)\times ...\times C^{\infty }(M)\longrightarrow
C^{\infty }(M^{A},A),(f_{1},...,f_{p})\longmapsto \lbrack \eta
(f_{1},...,f_{p})]^{A}  \notag
\end{equation}%
is skew-symmetric multilinear forms of degree $p$ from $C^{\infty }(M)$ into 
$C^{\infty }(M^{A},A)$ i.e 
\begin{equation*}
\eta ^{A}\in \Lambda _{Pois}^{p}(M^{A},\sim )\text{.}
\end{equation*}%
\newline

Thus

\begin{proposition}
For any $\eta\in\Lambda_{Pois}^{p}(M)$, we have $\widetilde{d}%
\eta^{A}=\left( d\eta\right) ^{A}$.
\end{proposition}

\begin{proof}
For any $f_{1},...,f_{p+1}\in C^{\infty}(M),$ we have 
\begin{align*}
(\widetilde{d}\eta^{A})(f_{1},...,f_{p+1}) & =\sum_{i=1}^{p+1}(-1)^{i}%
\widetilde{[ad(f_{i})]^{A}}\left( \eta^{A}(f_{1},...,\widehat{f_{i}}%
,...,f_{p+1})\right) \\
& +\sum_{1\leq i<j\leq p+1}(-1)^{i+j}\eta^{A}\left(
\{f_{i},f_{j}\},f_{1},...,\widehat{f_{i}},...,\widehat{f_{j}}%
,...,f_{p+1}\right) \\
& =\sum_{i=1}^{p+1}(-1)^{i}\widetilde{[ad(f_{i})]^{A}}\left( \eta (f_{1},...,%
\widehat{f_{i}},...,f_{p+1})\right) ^{A} \\
& +\sum_{1\leq i<j\leq p+1}(-1)^{i+j}[\eta(\{f_{i},f_{j}\},f_{1},...,%
\widehat{f_{i}},...,\widehat{f_{j}},...,f_{p+1})]^{A} \\
& =\sum_{i=1}^{p+1}(-1)^{i}\{f_{i},\eta(f_{1},...,\widehat{f_{i}}%
,...,f_{p+1})\}^{A} \\
& +\sum_{1\leq i<j\leq p+1}(-1)^{i+j}[\eta(\{f_{i},f_{j}\},f_{1},...,%
\widehat{f_{i}},...,\widehat{f_{j}},...,f_{p+1})]^{A} \\
& =[(d\eta)(f_{1},f_{2},...,f_{p+1})]^{A}\text{.}
\end{align*}

That ends the proof.
\end{proof}

\begin{corollary}
The $1$-form $\eta^{A}$ is $\widetilde{d}$-closed i.e ($\widetilde{d}\eta
^{A}=0$), if only if $d\eta=0$. In particular when $\eta$ is a derivation of
Poisson algebra $C^{\infty}(M)$.
\end{corollary}

\bigskip

\begin{proof}
Indeed, for $p=1$, we have%
\begin{align*}
(\widetilde{d}\eta^{A})(f,g) & =\widetilde{[ad(f)]^{A}}[\eta^{A}(g)]-%
\widetilde{[ad(g)]^{A}}[\eta^{A}(f)]-\eta^{A}(\{f,g\}) \\
& =\{f,\eta(g)\}^{A}-\{g,\eta(f)\}^{A}-[\eta({f,g})]^{A} \\
& =(\{f,\eta(g)\}-\{g,\eta(f)\}-[\eta({f,g})])^{A} \\
& =[d\eta(\{f,g)\}]^{A}\text{.}
\end{align*}
for any $f,g\in C^{\infty}(M)$.\bigskip

Thus $\widetilde{d}\eta^{A}=0$ if only if $d\eta=0$.

When $\eta$ is a derivation of Poisson algebra $C^{\infty}(M)$, we have $%
f,g\in C^{\infty}(M)$, 
\begin{align*}
\eta(\{f,g\}) & =\{f,\eta(g)\}-\{g,\eta(f)\} \\
& =\{\eta(f),g\}+\{f,\eta(g)\}
\end{align*}
i.e 
\begin{align*}
(\widetilde{d}\eta^{A})(f,g) & =[d\eta(\{f,g)\}]^{A} \\
& =0\text{.}
\end{align*}
That ends the proof.
\end{proof}

\begin{proposition}
If $\eta$ and $\eta^{\prime}$ both are cohomologous $d$-closed $p$-forms
then $\eta^{A}$ and $\eta^{\prime}{}^{A}$ both are cohomologous $\widetilde{d%
}$-closed $p$-forms.
\end{proposition}

\begin{proof}
For any $f_{1},...,f_{p}\in C^{\infty}(M)$ we have 
\begin{align*}
\lbrack\eta^{A}-\eta^{\prime}{}^{A}](f_{1},...,f_{p}) &
=\eta^{A}(f_{1},...,f_{p})-\eta^{\prime}{}^{A}(f_{1},...,f_{p}) \\
& =[\eta(f_{1},...,f_{p})]^{A}-[\eta^{\prime}(f_{1},...,f_{p})]^{A} \\
& =[\eta(f_{1},...,f_{p})-\eta^{\prime}(f_{1},...,f_{p})]^{A} \\
& =[(\eta-\eta^{\prime})(f_{1},...,f_{p})]^{A}\text{.}
\end{align*}

If there exists $\nu\in\Lambda_{Pois}^{p-1}(M)$ such that%
\begin{equation}
\eta-\eta^{\prime}=d\nu  \notag
\end{equation}
then%
\begin{align*}
\lbrack\eta^{A}-\eta^{\prime}{}^{A}](f_{1},...,f_{p}) & =[(\eta-\eta
^{\prime})(f_{1},...,f_{p})]^{A} \\
& =[d\nu(f_{1},...,f_{p})]^{A} \\
& =\widetilde{d}\nu^{A}(f_{1},...,f_{p}).
\end{align*}
i.e%
\begin{equation}
\eta^{A}-\eta^{\prime}{}^{A}=\widetilde{d}\nu^{A}\text{.}  \notag
\end{equation}
\end{proof}

The cohomology class of the $d$-closed $p$-form $\eta$ induces the
cohomology class of the $\widetilde{d}$-closed $p$-form $\eta^{A}$.

$\ $Let $Z_{Pois}^{p}(M^{A},\sim)$ be the set of $\widetilde{d}$-closed $p$%
-forms from $C^{\infty}(M)$ into $C^{\infty}(M^{A},A)$ and $%
B_{Pois}^{p}(M^{A},\sim)$ be the set of $\widetilde{d}$-exact $p$-forms from 
$C^{\infty}(M)$ into $C^{\infty}(M^{A},A)$. We denote 
\begin{equation*}
H_{Pois}^{p}(M^{A},\sim)=Z_{Pois}^{p}(M^{A},\sim)/B_{Pois}^{p}(M^{A},\sim)%
\text{.}
\end{equation*}

For $p=0,\Lambda_{Pois}^{0}(M^{A},\sim)=C^{\infty}(M^{A},A).$ It is obvious
that $H^{0}(M^{A},\sim)$ is the center of $C^{\infty}(M^{A},A)$ i.e the set 
\begin{equation*}
\left\{ \phi\in(M^{A},A);\widetilde{[ad(f)]^{A}}(\phi)=0\text{ for every }%
f\in C^{\infty}(M)\right\} \text{.}
\end{equation*}

For $p=1$,$\ $we have%
\begin{equation*}
H_{Pois}^{1}(M^{A},\sim)=0\text{.}
\end{equation*}

\section{$\widetilde{d_{A}}$-Poisson cohomology}

The map%
\begin{equation}
C^{\infty }(M^{A},A)\longrightarrow Der_{A}[C^{\infty }(M^{A},A)],\varphi
\longmapsto \widetilde{\tau _{\varphi }},  \notag
\end{equation}%
is a representation of $C^{\infty }(M^{A},A)$ into $C^{\infty }(M^{A},A).$
We denote $\widetilde{d_{A}}$ the cohomology operator associated to this
representation.\newline
For any $p\in \mathbb{N}$, $\Lambda _{Pois}^{p}(M^{A},\sim _{A})=\mathcal{C}%
^{p}[C^{\infty }(M^{A},A),C^{\infty }(M^{A},A)]$denotes the $C^{\infty
}(M^{A},A)$-module of skew-symmetric multilinear forms of degree $p$ on $%
C^{\infty }(M^{A},A)$ into $C^{\infty }(M^{A},A)$. We have%
\begin{equation}
\Lambda _{Pois}^{0}(M^{A},\sim _{A})=C^{\infty }(M^{A},A).  \notag
\end{equation}%
We denote%
\begin{equation}
\Lambda _{Pois}(M^{A},\sim _{A})=\bigoplus_{p=0}^{n}\Lambda
_{Pois}^{p}(M^{A},\sim _{A}).  \notag
\end{equation}%
For $\Omega \in \Lambda _{Pois}^{p}(M^{A},\sim _{A})$ and $\varphi
_{1},\varphi _{2},...,\varphi _{p+1}\in C^{\infty }(M^{A},A)$, we have%
\begin{align*}
\widetilde{d_{A}}\Omega (\varphi _{1},...,\varphi _{p+1})&
=\sum_{i=1}^{p+1}(-1)^{i-1}\widetilde{\tau _{\varphi _{i}}}[\Omega (\varphi
_{1},...,\widehat{\varphi _{i}},...,\varphi _{p+1}] \\
& +\sum_{1\leq i<j\leq p+1}(-1)^{i+j}\Omega (\widetilde{\tau _{\varphi _{i}}}%
(\varphi _{j}),\varphi _{1},...,\widehat{\varphi _{i}},...,\widehat{\varphi
_{j}},...,\varphi _{p+1})
\end{align*}%
i.e 
\begin{align*}
\widetilde{d_{A}}\Omega (\varphi _{1},\varphi _{2},...,\varphi _{p+1})&
=\sum_{i=1}^{p+1}(-1)^{i-1}\{\varphi _{i},\Omega (\varphi _{1},...,\widehat{%
\varphi _{i}},...,\varphi _{p+1}\}_{A} \\
& +\sum_{1\leq i<j\leq p+1}(-1)^{i+j}\Omega (\{\varphi _{i},\varphi
_{j}\}_{A},\varphi _{1},...,\widehat{\varphi _{i}},...,\widehat{\varphi _{j}}%
,...,\varphi _{p+1})\text{.}
\end{align*}%
For $p=1,$ we have 
\begin{equation*}
\widetilde{d_{A}}\Omega (\varphi ,\psi )=\{\Omega (\varphi ),\psi
\}_{A}+\{\varphi ,\Omega (\psi )\}_{A}-\Omega (\{\varphi ,\psi \}_{A})
\end{equation*}%
for any $\varphi ,\psi \in C^{\infty }(M^{A},A)$. Thus

\begin{corollary}
The $1$-form $\Omega$ is $\widetilde{d_{A}}$-closed i.e $\widetilde{d_{A}}%
\Omega=0$ if, and only if,%
\begin{equation*}
\Omega(\{\varphi,\psi\}_{A}=\{\Omega(\varphi),\psi\}_{A}+\{\varphi,\Omega
(\psi)\}_{A}
\end{equation*}
i.e $\Omega$ is a derivation of the algebra $C^{\infty}(M^{A},A)$.
\end{corollary}

\bigskip$\ $Let $Z_{Pois}^{p}(M^{A},\sim_{A})$ be the set of $\widetilde {%
d_{A}}$-closed $p$-forms from $C^{\infty}(M^{A},A)$ into $%
C^{\infty}(M^{A},A) $ and $B_{Pois}^{p}(M^{A},\sim_{A})$ be the set of $%
\widetilde{d_{A}}$-exact $p$-forms from $C^{\infty}(M)$ into $%
C^{\infty}(M^{A},A)$. We denote 
\begin{equation*}
H_{Pois}^{p}(M^{A},\sim_{A})=Z_{Pois}^{p}(M^{A},\sim_{A})\
/B_{Pois}^{p}(M^{A},\sim_{A})\text{.}
\end{equation*}

For $p=0,\Lambda _{Pois}^{0}(M^{A},\sim _{A})=C^{\infty }(M^{A},A).$ It is
obvious that $H^{0}(M^{A},\sim _{A})$ is the center of $C^{\infty }(M^{A},A)$
i.e the set 
\begin{equation*}
\left\{ \varphi \in C^{\infty }(M^{A},A);\left\{ \varphi ,\phi \right\}
_{A}=0\text{ for every }\phi \in C^{\infty }(M^{A},A)\right\} \text{.}
\end{equation*}

For $p=1$,$\ $we have 
\begin{equation*}
H_{Pois}^{1}(M^{A},\sim _{A})=0\text{.}
\end{equation*}

\bigskip

\textbf{Acknowledgements: } The first author thanks Deutscher Akademischer
Austauschdientst (DAAD) for their financial support.


\bigskip

\begin{thebibliography}{99}
\bibitem{bok} B. G. R. Bossoto, E. Okassa, \textit{Champs de vecteurs et
formes diff\'{e}rentielles sur une vari\'{e}t\'{e} de points proches},
Archivum Mathematicum (Brno), 44(2008) 159-171.

\bibitem{bos2} B. G. R. Bossoto, \textit{Structures de Jacobi sur une vari%
\'{e}t\'{e} des points proches}, Math. Vesnik. 62, 2 (2010), 155-167.

\bibitem{bok1} B. G. R. Bossoto, E. Okassa, \textit{A-poisson structures on
Weil bundles}, Int., J. Contemp. Math. Sciences, Vol. 7, 2012, n%
${{}^\circ}$%
16, 785-803.

\bibitem{kms} I. Kolar, P.W. Michor, and J. Slovak, \textit{Natural
operations in differential geometry}. Springer, 1993, 434 p.

\bibitem{Lic} A. Lichnerowicz, \textit{Les vari\'{e}t\'{e}s de Poisson et
leurs alg\'{e}bres de Lie associ\'{e}es}, J. Diff. Geom., 12 (1977),
253--300.

\bibitem{mor} A. Morimoto, \textit{Prolongation of connections to bundles of
infinitely near points}, J. Diff. Geom, t.11(1976), 479-498.

\bibitem{oka} E. Okassa, \textit{Prolongement des champs de vecteurs \`{a}
des vari\'{e}t\'{e}s de points proches}, Annales Facult\'{e} des sciences de
Toulouse, Vol. VIII, n$^{\circ }$3, 1986-1987, 349-366.

\bibitem{iva} I. Vaisman, \textit{Lectures on the Geometry of Poisson
Manifolds}, Progress in Math.118, Birkh\"{a}user Verlag, Basel, 1994.

\bibitem{shu} V. V. Shurygin, \textit{Smooth manifolds over local algebras
and Weil bundles}, J. Math. Sci., 108 (2) (2002), 249--294.

\bibitem{wei} A. Weil, \textit{Th\'{e}orie des points proches sur les vari%
\'{e}t\'{e}s diff\'{e}rentiables}, Colloq. G\'{e}om. Diff. Strasbourg
(1953), 111-117.
\end{thebibliography}
\end{document}